\documentclass[leqno,draft,11pt]{article}
\usepackage{amssymb,amsmath,theorem,a4wide}

\def\shortthm{}
\def\enumup{}

\def\dedic{\thanks{Supported by Institutional Research Plan
AV0Z10190503, grant IAA100190902 of GA AV \v CR, and project 1M0545 of
M\v SMT \v CR.}}

\ifx\mydefs\undefined\else \fi
\let\mydefs\relax

\allowdisplaybreaks[2]

\ifx\loadcyr\undefined\else
\input cyracc.def
\makeatletter
 at 1\@ptsize pt
 at 1\@ptsize pt
 at 1\@ptsize pt
\makeatother

\fi

\let\M\mathit
\def\gobble#1{}
\def\fixsup#1#2{{#1\let\dp\gobble\mathstrut}^#2_}

\def\bme{\hskip.75em\relax}

\let\eq\leftrightarrow

\def\?{\mathbin?}

\let\model\vDash

\newbox\circlebox
\setbox\circlebox\hbox{$\bigcirc$}
\def\circled#1{%
  \setbox0\hbox to\wd\circlebox{\hss$#1$\hss}\wd0=0pt
  \box0\copy\circlebox}

\let\tet\vartheta
\let\ep\varepsilon

\def\greek#1{$\expandafter\greeknum\csname c@#1\endcsname$}
\makeatletter
\def\greeknum#1{\ifcase#1\or\alpha\or\beta\or\gamma\or\delta\or\ep
      \or\digamma\or\zeta\or\eta\or\tet\or\iota\else\@ctrerr\fi}
\makeatother

\def\p#1{\langle#1\rangle}

\let\onto\twoheadrightarrow

\def\cupd{\mathbin{\dot\cup}}

\def\twoprimes{\raise.2\fontdimen6\scriptfont2\hbox{$\scriptstyle\prime\prime$}}
\newcommand\rpair[3][3em]{\mathrel{%
   \begin{matrix}%
     \strut\smash{\xrightonto{\hbox to#1{\hss$#2$\hss}}}\\[-1.7ex]%
     \strut\smash{\xleftembed[\hbox to#1{\hss$#3$\hss}]{}}%
   \end{matrix}}}
\makeatletter
\newcommand\xrightonto[2][]{\ext@arrow 0359\rightontofill{#1}{#2}}
\newcommand\xleftembed[2][]{\ext@arrow 3095\leftembedfill{#1}{#2}}
\def\leftembedfill{\arrowfill@\leftarrow\relbar\hookleftnoarrow}
\def\rightontofill{\arrowfill@\relbar\relbar\onto}
\def\hookleftnoarrow{\DOTSB\relbar\joinrel\rhook}
\makeatother

\def\rem{\mathbin{\mathrm{rem}}}

\mathchardef\#="2023 

\ifx\busspf\undefined\else
\usepackage{bussproofs}
\EnableBpAbbreviations

\fi

\ifx\symlasy\undefined

  \def\centerdot#1{%
    \setbox0\hbox{$\mathop{#1}$}\dimen0 \ht0
    \setbox0\hbox{$#1$}\advance\dimen0 -\ht0
    \setbox2\hbox to\wd0{\hss$\mathop{\cdot}$\hss}\wd2=0pt
    \lower\dimen0\box2\box0 }
\else
  
  \def\centerdot#1{%
     \setbox0\hbox{$#1$}%
     \raise0.206\ht0\hbox to\wd0{\hss$\cdot$\hss}%
     \kern-\wd0 \box0 }
\fi

\let\sls|

\def\Up{{\setbox0\hbox{$\uparrow$}%
         \lower\dp0\hbox to\wd0{\hss\vrule width4pt height.4pt\hss}%
         \kern-\wd0\box0}}
\def\UP{{\setbox0\hbox{$\uparrow$}%
         \lower\dp0\hbox to\wd0{\hss\vrule width4pt height.4pt\hss}%
         \kern-\wd0\copy0\kern-\wd0\raise.35ex\box0}}
\def\Down{{\setbox0\hbox{$\downarrow$}%
         \raise\ht0\hbox to\wd0{\hss\vrule width4pt depth.4pt\hss}%
         \kern-\wd0\box0}}

\def\st{\expandafter\hat}

\def\PA{\M{PA}}

\let\stm\N

\mathcode`\*="0203

\mathchardef\mhyphen="2D

\def\noproof{\leavevmode\unskip\bme\vadjust{}\nobreak\hfill$\qed$\par}
\let\qed\Box
\newenvironment{Pf}[1][]
  {\par\noindent\textit{Proof\optpar{#1}:}\bme\ignorespaces}
  {\noproof\pagebreak[2]\vskip\medskipamount\ignorespacesafterend}
\def\qedhere{\relax\ifmmode\eqno\qed\expandafter\aftergroup
                   \else\noproof\fi\noqed}
\def\noqed{\let\noproof\relax}

\theoremstyle{plain}
\ifx\shortthm\undefined
\newtheorem{Thm}{Theorem}[section]
\else
\newtheorem{Thm}{Theorem}
\fi

\newtheorem{Lem}[Thm]{Lemma}

\def\theCl{\arabic{Cl}}

\theorembodyfont\upshape
\newtheorem{Def}[Thm]{Definition}

\newenvironment{Pf*}{\let\qed\qedCl\Pf}\endPf

\usepackage[reftex]{theoremref}

{\catcode`\^^I=13 \catcode`\^^M=13
\gdef\doalgo#1#2\end#{\hbox to\hsize{\hss \let^^I\qquad%
  \def\\^^M{\nobreak\hfil\break\vadjust{}\qquad}%
  \fboxsep1em \linenum0 %
  \fbox{\hsize#1\vbox{%
  \everypar{\advance\linenum1 %
      \hbox to1em{\hss$\scriptstyle\the\linenum$}\quad}%
  #2}}\hss}\end}}
\newcount\linenum

\def\key{\relax\ifmmode\expandafter\mathbf\else\expandafter\textbf\fi}

\def\allowhyphens{\nobreak\hskip0pt\relax}

\DeclareRobustCommand*\magiclparen{\ifmmode(\else\textup(\allowhyphens\fi}
\DeclareRobustCommand*\magicrparen{\ifmmode)\else\textup)\fi}
\let\lparen=(  \let\rparen=)
\def\magicparon{\catcode`\(\active\catcode`\)\active}
\def\magicparoff{\catcode`\(12 \catcode`\)12 }
\AtBeginDocument{\ifx\ifPreview\iftrue\else\magicparon\fi}
\magicparon
\let (=\magiclparen  \let )=\magicrparen

\ifx\enumup\undefined
  
\else
  
\fi

\def\optpar#1{\ifx\relax#1\relax\else\/ (#1)\fi}

\magicparoff

\mathchardef\comma=\mathcode`\,
{\catcode`\,=\active \gdef,{\comma\penalty\relpenalty}}

\providecommand\dedic{}
\author{Emil Je\v r\'abek\dedic\\[\medskipamount]
Institute of Mathematics of the Academy of Sciences\\
\small \v Zitn\'a 25,
115\:67 Praha 1,
Czech Republic,
email: \texttt{jerabek@math.cas.cz}
}

\title{Sequence encoding without induction}

\begin{document}
\maketitle
\begin{abstract}
We show that the universally axiomatized, induction-free theory
$\PA^-$ is a sequential theory in the sense of Pudl\'ak~\cite{pud:cuts}, in
contrast to the closely related Robinson's arithmetic.
\end{abstract}
Ever since G\"odel's \cite{godel} arithmetization of syntax in the proof of his
incompleteness theorem, sequence encoding has been an indispensable
tool in the study of arithmetical theories and related areas of
mathematical logic. While a common approach is to develop a particular
sequence encoding in a suitable base theory and work with that, a
general concept of theories supporting encoding of sequences of their
elements, called \emph{sequential theories}, was isolated by
Pudl\'ak~\cite{pud:lattice,pud:cuts} during his work on
interpretability.
A similar but weaker notion was defined earlier by Vaught~\cite{vaught}.
More generally, theories with ``containers'' of various kind
were studied by Visser~\cite{vis:pair}, who includes a discussion of
variants of the notion of sequentiality and some historical remarks.

It is known that fairly weak arithmetical theories can be sequential
(e.g., fragments of bounded arithmetic such as Buss' $S^1_2$, cf.\
Kraj\'\i\v cek~\cite{book}), nevertheless sequential theories
described in the literature so far generally involve some form of the induction
schema. Indeed, the common induction-free base arithmetical theory,
Robinson's $Q$ \cite{tmr}, is \emph{not} sequential (Visser~\cite{vis:pair}).
This also shows that sequentiality is not preserved by interpretations
(even in the tame form of definable cuts): the sequential theory
$S^1_2$ is interpretable on a cut in $Q$. The reason is that in a
sequential theory, all elements in the universe of the theory have to
be admissible as sequence entries, not just elements from a proper
cut. (On the other hand, the \emph{lengths} of sequences may be
confined to a small cut.)

In this note, we prove sequentiality of the theory $\PA^-$ of
discretely ordered commutative semirings with the least element
(without the subtraction axiom), and therefore of all its simple
extensions. Like $Q$, $\PA^-$ is an 
induction-free theory axiomatized by basic properties of
$+,\cdot,\le$, and it (or its slightly stronger variants) is often
used as an arithmetical base theory (see e.g.\
Kaye~\cite{kaye:models}). Our version of $\PA^-$ is purely universally
axiomatized. The main result can also be adapted in a straightforward way
to the theory of discretely ordered (commutative) rings.

We encode sequences in $\PA^-$ using the well-known G\"odel's
$\beta$-function, slightly modified for a technical reason. Where the
usual analysis of G\"odel's $\beta$ employs induction, we switch to a
shorter cut; the main problem is to ensure we can make do with
restricting only the lengths of sequences to the cut, while allowing
arbitrary elements to appear in sequences.
We proceed with the formal details.

\begin{Def}\th\label{def:p-}
Let $\PA^-$ be the theory of discretely ordered commutative semirings
with the least element. That is, $\PA^-$ is the first-order theory
with equality in the language
$\p{0,1,{+},{\cdot},{\le}}$, axiomatized by
\begin{gather}
\tag{A1}x+0=x\\
\tag{A2}x+y=y+x\\
\tag{A3}(x+y)+z=x+(y+z)\\
\tag{M1}x\cdot1=x\\
\tag{M2}x\cdot y=y\cdot x\\
\tag{M3}(x\cdot y)\cdot z=x\cdot(y\cdot z)\\
\tag{AM}x\cdot(y+z)=x\cdot y+x\cdot z\\
\tag{O1}x\le y\lor y\le x\\
\tag{O2}(x\le y\land y\le z)\to x\le z\\
\tag{S1}x+1\nleq x\\
\tag{S2}x\le y\to(x=y\lor x+1\le y)\\
\tag{OA}x\le y\to x+z\le y+z\\
\tag{OM}x\le y\to x\cdot z\le y\cdot z
\end{gather}
Let $x<y$ abbreviate $x\le y\land x\ne y$.

Note that many authors (e.g., Kaye~\cite{kaye:models} or
Kraj\'\i\v cek~\cite{book}) use
a stronger definition of $\PA^-$, namely as the theory of nonnegative parts of
discretely ordered rings, which includes the subtraction axiom $x\le
y\to\exists z\,(z+x=y)$. In contrast, our version of
$\PA^-$ is a universal theory, hence it does not even prove the
existence of predecessors (e.g., the semiring $\stm[x]$ of polynomials
with nonnegative integer coefficients, ordered lexicographically, is a
model of $\PA^-$).

Sequentiality can be defined in several ways. For definiteness, we will follow
the (relatively restrictive) definition of Pudl\'ak~\cite{pud:cuts}: a
theory $T$ is \emph{sequential} if it contains Robinson's arithmetic
$Q$ relativized to some formula $N(x)$, and there is a formula
$\beta(x,i,w)$ (whose intended meaning is that $x$ is the $i$th
element of a sequence $w$) such that $T$ proves
\[\tag{SEQ}
\forall w,x,k\,\exists w'\,\forall i,y\,[(N(k)\land i\le k)\to
  [\,\beta(y,i,w')\eq((i<k\land\beta(y,i,w))\lor(i=k\land y=x))]].\]

A definable set is called \emph{inductive} if it contains $0$ and is
closed under successor, and it is a \emph{cut} if it is furthermore
downward closed.
\end{Def}

We first establish some basic properties of $\PA^-$ which the reader
might have been missing among the axioms:
\begin{Lem}\th\label{lem:p-}
$\PA^-$ proves
\begin{enumerate}
\item $(x\le y\land y\le x)\to x=y$
\item $x+z\le y+z\to x\le y$
\item $x\cdot0=0$
\item $0\le x$
\item $(z\ne 0\land x\cdot z\le y\cdot z)\to x\le y$
\item $x\le y+1\eq(x\le y\lor x=y+1)$
\end{enumerate}
\end{Lem}
\begin{Pf}
(i): Otherwise $x+1\le y$ by S2, hence $x+1\le x$ by O2, contradicting
S1.

(ii): Otherwise $y<x$ by O1, hence $y+1\le x$ by S2 and
$(x+z)+1\le(y+z)+1=(y+1)+z\le x+z$ by OA, A2, A3, contradicting (O2
and) S1.

(iii): $x\cdot0+x=x\cdot0+x\cdot1=x(0+1)=x\cdot1=x=0+x$, hence $x\cdot0=0$ by (ii) and (i).

(iv): $0\le 1$ by S1 and O1, hence $0=x\cdot0\le x\cdot1=x$.

(v): Otherwise $y<x$, thus $y+1\le x$, giving $yz\le yz+z=(y+1)z\le
xz\le yz$, hence $yz=yz+z$ by (i), and $z=0$ by (ii) and (i).

(vi): Left-to-right: if $x<y+1$, then $x+1\le y+1$ by S2, hence $x\le y$.
\end{Pf}

\begin{Def}\label{def:beta}
Let $\p{x,y}:=(x+y)^2+x$, and let
$$\beta(x,i,w):\eq\exists u,v,q\,[w=\p{u,v}\land
u=q(1+(i+1)v)+x\land x\le(i+1)v]$$
be G\"odel's $\beta$-function.
\end{Def}

\begin{Lem}\th\label{lem:uniq}
$\PA^-$ proves:
\begin{enumerate}
\item $(dx+y=dx'+y'\land y,y'<d)\to x=x'\land y=y'$
\item $\p{x,y}=\p{x',y'}\to(x=x'\land y=y')$
\item $(\beta(x,i,w)\land\beta(x',i,w))\to x=x'$
\end{enumerate}
\end{Lem}
\begin{Pf}
(i): If $x<x'$, then $dx+y<dx+d=d(x+1)\le dx'\le dx'+y'=dx+y$, a
contradiction. Thus $x\ge x'$, and symmetrically, $x'\ge x$, thus
$x=x'$, which implies $y=y'$.

(ii): We have $(x+y)^2\le\p{x,y}\le(x+y)^2+(x+y)<(x+y+1)^2$, and $u^2$
is monotone. Thus, $\p{x,y}=\p{x',y'}$ implies $x+y=x'+y'$, which in
turn implies $x=x'$, which implies $y=y'$.

(iii): $u,v$ are unique by (ii), and then ($q$ and) $x$ is unique by (i).
\end{Pf}

\begin{Def}\th\label{def:mod}
Let $x\rem y=z$ denote $z<y\land
\exists q\,(x=z+qy)$ (this is $\PA^-$-provably a partial function by
\th\ref{lem:uniq}). We write $y\mid x$ for $\exists q\,(x=qy)$. 
\end{Def}

Note that $\beta(x,i,\p{u,v})$ iff $x=u\rem(1+(i+1)v)$.

One step in the usual proof that G\"odel's $\beta$-function works
is to show that the numbers $1+v$, $1+2v$, \dots, $1+kv$ are pairwise
coprime if $v$ is divisible by $1$, \dots, $k-1$. The next lemma can be
vaguely thought of as a replacement for this statement in our situation.
\begin{Lem}\th\label{lem:i1}
Let
\begin{align*}
I_0(k)&:\eq\forall j\le i\le k\,\exists d\,(d+j=i),\\
I_1(k)&:\eq I_0(k)\land\forall v\,\exists u\,
  [\,\forall i\le k\,(1+iv\mid u)\\
&\qquad\qquad\quad\land\forall i>k\,[(I_0(i)\land\forall0<j\le k\,(i-j\mid v))
      \to\exists p\,(up\rem(1+iv)=1)]].
\end{align*}
Then $\PA^-$ proves that $I_0$ is a cut and $I_1$ is inductive. Here,
$i-j$ denotes the (unique) $d$ such that $d+j=i$, which exists because
of $I_0(i)$.
\end{Lem}
\begin{Pf}
That $I_0$ is a cut is easy to see. $I_1(0)$ follows by taking $u=1$.
Assume $I_1(k)$, and let $v$ be given. Let $u$ be the witness for
$I_1(k)$, and put $u'=(1+(k+1)v)u$. Clearly $1+iv\mid u'$ for all
$i\le k+1$. Let $i>k+1$ be such that $I_0(i)$ and $i-j\mid v$ for all
$j\le k+1$, $j>0$. By $I_1(k)$, there exist $p,q$ such that
$up=1+(1+iv)q$.
Moreover, we claim that
\[\tag{$*$}\bigl(1+(k+1)v\bigr)p'=1+(1+iv)q'\]
for some $p',q'$. Then $u'pp'=1+(1+iv)\bigl(q+q'+qq'(1+iv)\bigr)$,
which completes the proof of $I_1(k+1)$.

In order to show $(*)$, write $k'=k+1$, and fix $z$ such that
$(i-k')z=v$, which exists by our assumption on $i$. We have
$k'v+k'^2z=ik'z$, hence
$$(1+k'v)+(1+iv)k'^2z=1+(1+k'v)ik'z.$$
We add a suitable multiple of $(1+k'v)(1+iv)$ to both sides in order
to move $1+iv$ to the right-hand side and $1+k'v$ to the left-hand
side, as required in $(*)$:
\begin{align*}
(1+k'v)\bigl(1&+k'+ik'(i-(k'+1))z\bigr)
  +(1+k'v)ik'z+(1+iv)k'^2z\\
&=(1+k'v)(1+k'+ik'v)+(1+iv)k'^2z\\
&=(1+k'v)+(1+iv)k'^2z+(1+k'v)(1+iv)k'\\
&=1+(1+k'v)ik'z+(1+iv)(1+k'v)k'\\
&=1+(1+iv)\bigl(k'+k'^2(i-(k'+1))z\bigr)+(1+iv)k'^2z+(1+k'v)ik'z,
\end{align*}
and we can cancel $(1+k'v)ik'z+(1+iv)k'^2z$ from both sides using
\th\ref{lem:p-}. Thus, we have $(*)$ with $p'=1+k'+ik'(i-(k'+1))z$ and
$q'=k'+k'^2(i-(k'+1))z$.
\end{Pf}
The main point of the following definition of $I_2(k)$ is that
$\beta$-encoded sequences of length $k$ can be recoded using a
different $v$, as well as expanded by a $(k+1)$th element. This is
made more explicit in \th\ref{lem:beta}.
\begin{Lem}\th\label{lem:i2}
Define
\begin{align*}
I_2(k):\eq I_1(k)&\land\forall u,v,v',x\,\exists u'\\
  [(\forall0<i&\le k\,(i\mid v')\land v'\ge v\land(k+1)v'\ge x\land\forall i\le k\,
             \exists r\,(r=u\rem(1+iv)))\\
   \to(&\forall i\le k\,\exists r\,
           (r=u\rem(1+iv)\land r=u'\rem(1+iv'))\\
      &\land x=u'\rem(1+(k+1)v'))].
\end{align*}
Then $\PA^-$ proves that $I_2$ is inductive.
\end{Lem}
\begin{Pf}
For $k=0$, we can take $u':=x$.

Assume $I_2(k)$, we will prove $I_2(k+1)$. Let $u,v,v',x$ be given. By
$I_2(k)$, we can find a $u_0$ such that there exists
$u_0\rem(1+iv')=u\rem(1+iv)$ for all $i\le k+1$, using the fact that
$u\rem(1+(k+1)v)\le(k+1)v\le(k+1)v'$. Since $I_1(k+1)$ and $v'$ is
divisible by $1,\dots,k+1$, there are $u_1,p,q$ such that $1+iv'\mid
u_1$ for all $i\le k+1$, and $u_1p=1+(1+(k+2)v')q$. Define
$u':=u_0+(x+u_0(k+2)v')u_1p$. Then $u'\rem(1+iv')=u\rem(1+iv)$ for all
$i\le k+1$, and $x=u'\rem(1+(k+2)v')$, as
\begin{align*}
u'&=u_0+\bigl(x+u_0(k+2)v'\bigr)\bigl(1+(1+(k+2)v')q\bigr)\\
&=\bigl(x+u_0(k+2)v'\bigr)q\bigl(1+(k+2)v'\bigr)+u_0+u_0(k+2)v'+x\\
&=\bigl(u_0+xq+u_0(k+2)v'q\bigr)\bigl(1+(k+2)v'\bigr)+x,
\end{align*}
and $x<1+(k+2)v'$.
\end{Pf}

\begin{Lem}\th\label{lem:beta}
$\PA^-$ proves: if $I_2(k)$ and $\forall i<k\,\exists x\,\beta(x,i,w)$,
then there exists a $w'$ such that
$$\forall i\le k\,\forall y\,
  [\,\beta(y,i,w')\eq((i<k\land\beta(y,i,w))\lor(i=k\land y=x))].$$
\end{Lem}
\begin{Pf}
Let $w=\p{u_0,v_0}$, and write $(w)_i=x$ instead of $\beta(x,i,w)$ for
clarity. Applying $I_1(k)$ with $v=1$, we see that there exists a
$v'>0$ divisible by $1,\dots,k+1$. Pick $v_1$ such that $v_1\ge v_0$,
$v_1\ge x$, and
$v'\mid v_1$. By $I_2(k)$, there exists $u_1$ such that for all $i<k$,
$u_1\rem(1+(i+1)v_1)=u_0\rem(1+(i+1)v_0)$, and $u_1\rem(1+(k+1)v_1)=x$. Thus,
if we put $w':=\p{u_1,v_1}$, then $(w')_i=(w)_i$ for all $i<k$, and
$(w')_k=x$.
\end{Pf}
Clearly, \th\ref{lem:i2,lem:beta} almost show that $\PA^-$ is sequential.
However, as $\PA^-$ does not prove that division with remainder is
total, it may happen 
for G\"odel's $\beta$-function that $(w)_i$ is
undefined for some values of $i<k$, and then the definition of
sequentiality requires $(w')_i$ to be also undefined for the same
values of $i$. This does not seem possible to arrange, as we have no
way of forcing $(w')_i$ to be undefined when building $w'$. We fix this problem
by modifying the definition of $\beta$ a little bit.
\begin{Def}\th\label{def:betap}
$$\beta'(x,i,w):\eq[\,\beta(x,i,w)\land\forall j<i\,\exists y\,\beta(y,j,w)]
      \lor[x=0\land\exists j\le i\,\neg\exists y\,\beta(y,j,w)].$$
Note that $\beta'$ is $\PA^-$-provably a total function.
\end{Def}
\begin{Lem}\th\label{lem:i3}
$\PA^-$ proves that
$$I_3(k):\eq I_2(k)\land\forall w'\,\exists w\,\forall i<k\,\forall x\,
   (\beta(x,i,w)\eq\beta'(x,i,w'))$$
is inductive.
\end{Lem}
\begin{Pf}
$I_3(0)$ is clear. Assuming $I_3(k)$, we have $I_2(k+1)$ by
\th\ref{lem:i2}. Let $w$ be given, and write $x=(w)_i$ instead of
$\beta'(x,i,w)$ for clarity.
Since $I_3(k)$, there exists $w'$ such that $\beta((w)_i,i,w')$ for all
$i<k$. By \th\ref{lem:beta}, there exists $w''$ such that
$\beta((w)_i,i,w'')$ for $i<k$, and $\beta((w)_k,k,w'')$. This shows
$I_3(k+1)$.
\end{Pf}
By the usual shortening of cuts, let $N(x)$ be such that $\PA^-$ proves
that $N$ is a cut closed under $+$ and $\cdot$, and $N(x)\to I_3(x)$.
\begin{Thm}\th\label{thm:seq}
$\PA^-$ is a sequential theory with respect to $N$ and $\beta'$.
\end{Thm}
\begin{Pf}
$\PA^-$ proves itself relativized to $N$, as it is a universal theory.
Moreover, if $x\le y$ and $N(y)$, there exists $z$ such that
$z+x=y$ as $I_0(y)$, and we have $N(z)$ as $N$ is downward closed.
Thus, the subtraction axiom holds in $N$, hence $N$ is an
interpretation of $Q$ in $\PA^-$.

In order to show (SEQ) for $\beta'$, let $w,x,k$ such that $N(k)$ be given, and
write $(w)_i=y$ for $\beta'(y,i,w)$. By the definition of $I_3$, we
can find a $w'$ such that $\beta((w)_i,i,w')$ for each $i<k$. Then
\th\ref{lem:beta} gives a $w''$ such that $\beta((w)_i,i,w'')$ for
each $i<k$, and $\beta(x,k,w'')$. By the definition of $\beta'$, this
implies $(w'')_i=(w)_i$ for $i<k$, and $(w'')_k=x$.
\end{Pf}

In contrast, Robinson's $Q$ is not sequential, despite that it is
fairly close to $\PA^-$ in strength. In fact, Visser~\cite{vis:pair}
proved that it does not even support pairing; we include a
somewhat different proof of his result below for completeness:
\begin{Thm}\th\label{thm:q}
$Q$ is not sequential, and it has no pairing operation:
i.e., there is no formula $\pi(x,y,p)$ such that $Q$ proves
\begin{enumerate}
\item $\forall x,y\,\exists p\,\pi(x,y,p)$,
\item $\forall x,y,x',y',p\,[(\pi(x,y,p)\land\pi(x',y',p))\to(x=x'\land y=y')]$.
\end{enumerate}
\end{Thm}
\begin{Pf}
Let $M=\stm\cupd\{a_0,a_1\}$, and define arithmetical operations on
$M$ extending the usual operations on $\stm$ by $a_i+x=a_i$,
$n+a_i=a_i$, $n\cdot a_i=a_i$, $a_i\cdot0=0$, $a_i\cdot x=a_i$ if
$x\ne0$, where $x\in M$, and $n\in\stm$. Then $M\model Q$, and the
function $f$ identical on $\stm$ such that $f(a_i)=a_{1-i}$ is an
automorphism of $M$. Let $\pi$ be a pairing operation, and find an $x$
such that $\pi(a_i,a_j,x)$. Since $f$ is an automorphism,
$\pi(a_{1-i},a_{1-j},f(x))$. By unique decoding of pairs, it follows
that $f(x)\ne x$, i.e., $x\in\{a_0,a_1\}$. However, there are only two
elements in $\{a_0,a_1\}$, while there are four pairs of the form
$\p{a_i,a_j}$, contradicting uniqueness.
\end{Pf}
\subsection*{Acknowledgement}
I am grateful to Albert Visser for many useful suggestions.

\bibliographystyle{mybib}
\bibliography{mybib}

\providecommand\gobble[1]{} {\catcode`\/=13
  \gdef/{\string/\futurelet\nexttoken\finishslash}
  \gdef\finishslash{\ifx\nexttoken/\else\penalty\relpenalty\fi}}
  \providecommand\url{\begingroup\catcode`\~=12 \catcode`\/=13 \finishurl}
  \def\finishurl#1{\texttt{#1}\endgroup}
  \providecommand\dotminus{\mathbin{\setbox0\hbox{$-$}\setbox2\hbox
  to\wd0{\hss$^{\mkern1mu\cdot}$\hss}\wd2=0pt\box2\box0}}
\providecommand{\bysame}{\leavevmode\hbox to5em{\hrulefill}\thinspace}
\providecommand\bibliographyhook{}
\begin{thebibliography}{1}
\bibliographyhook

\bibitem{godel}
Kurt G{\"o}del, \emph{{\"U}ber for\-mal un\-ent\-scheid\-ba\-re {S}{\"a}tze der
  {P}rincipia {M}athematica und ver\-wand\-ter {S}y\-ste\-me, {I}},
  Mo\-nats\-hef\-te f{\"u}r {M}a\-the\-ma\-tik und {P}hy\-sik 38 (1931),
  pp.~173--198.

\bibitem{kaye:models}
Richard Kaye, \emph{Models of {P}eano arithmetic}, Oxford Logic Guides vol.~15,
  Oxford University Press, 1991.

\bibitem{book}
Jan Kraj{\'\i}{\v c}ek, \emph{Bounded arithmetic, propositional logic, and
  complexity theory}, Encyclopedia of Mathematics and Its Applications vol.~60,
  Cambridge University Press, 1995.

\bibitem{pud:lattice}
Pavel Pudl{\'a}k, \emph{Some prime elements in the lattice of interpretability
  types}, Transactions of the American Mathematical Society 280 (1983), no.~1,
  pp.~255--275.

\bibitem{pud:cuts}
\bysame, \emph{Cuts, consistency statements and interpretations}, Journal of
  Symbolic Logic 50 (1985), no.~2, pp.~423--441.

\bibitem{tmr}
Alfred Tarski, Andrzej Mostowski, and Rafael~M. Robinson, \emph{Undecidable
  theories}, North-Holland, Amsterdam, 1953.

\bibitem{vaught}
Robert~L. Vaught, \emph{Axiomatizability by a schema}, Journal of Symbolic
  Logic 32 (1967), no.~4, pp.~473--479.

\bibitem{vis:pair}
Albert Visser, \emph{Pairs, sets and sequences in first-order theories},
  Archive for Mathematical Logic 47 (2008), no.~4, pp.~299--326.

\end{thebibliography}
\end{document}